\documentclass{elsarticle}
\usepackage{hyperref}
\usepackage{mathrsfs}
\let\memoldbibsection\bibsection
\let\bibsection\relax
\usepackage{amsrefs}         
\let\bibsection\memoldbibsection
\usepackage{amsmath}

\def\int{\mathop{\rm int}\nolimits}

\newcommand{\N}{\mathbb{N}}
\newcommand{\Q}{\mathbb{Q}}
\newcommand{\U}{\mathcal{U}}
\makeatletter
\def\ps@pprintTitle{%
  \let\@oddhead\@empty
  \let\@evenhead\@empty
  \let\@oddfoot\@empty
  \let\@evenfoot\@oddfoot
}
\makeatother

\usepackage{amsfonts}
\usepackage{amssymb}
\newtheorem{thm}{Theorem}
\newtheorem{lem}[thm]{Lemma}
\newtheorem{exa}{Example}
\newtheorem{pro}{Proposition}
\newdefinition{rmk}{Remark}

\newtheorem{question}{Question}
\newtheorem{definition}{Definition}
\newproof{pf}{Proof}
\begin{document}
\begin{frontmatter}
\title{$Cp(X)$ for Hattori Spaces} 
\author[1]{Elmer Enrique Tovar Acosta\fnref{fn1,fn2}}
\ead{elmer@ciencias.unam.mx}

\affiliation[1]{organization={Universidad Nacional Autónoma de México.},
            addressline={Av. Universidad 3000}, 
            city={CDMX},
            postcode={4510}, 
            country={México}}  

\fntext[fn1]{The  author acknowledges financial support from CONACyT, grant no. 829699}
\fntext[fn2]{This paper is a part of the doctoral research of the author.}
\begin{abstract}
Motivated by the main results of the articles by Hattori \cite{Hattori} and Bouziad \cite{completez}, we seek to answer the following questions about Hattori spaces. Let $A\subseteq \mathbb{R}$, then:
\begin{enumerate}

    \item Given a compact set $K$ in the Euclidean topology, under what conditions is $K$ compact in the Hattori space $H(A)$?
    \item When is  $H(A)$ a quasi-metrizable space?
    \item When is  $H(A)$ a semi-stratifiable space?

\item When is $C_p(H(A))$ a normal space?
\item When is $C_p(H(A))$ a Lindelöf space?
\end{enumerate}
We obtain complete answers for 3 out of these 5 questions, while the last ones remain with partial answers, among them:\\ 
\

\noindent \textbf{Theorem}. If $\mathbb{R}\setminus A$ is analytic, then $C_p(\mathbb{R}_A)$ is not normal.\\

Moreover when we work on the Solovay Model we can improve the previous result to only require $\mathbb{R}\setminus A$ to be uncountable.
\end{abstract}
\end{frontmatter}
\section{Introduction}

In their article \cite{Hattori}, Hattori defines a family of intermediate topologies between the Euclidean and Sorgenfrey topologies using local bases. Namely, for every $A\subseteq \mathbb{R}$ he defines a topology $\tau(A)$ in $\mathbb{R}$ as follows: \begin{itemize}
    \item For every $x\in A$, the set $\lbrace (x-\varepsilon,x+\varepsilon) \ | \varepsilon>0\rbrace $ remains a local base at $x$.
    \item On the other hand, if $x\in \mathbb{R}\setminus A$, then $\lbrace [x,x+\varepsilon) \ | \ \varepsilon>0\rbrace$ is a local base at $x$.
\end{itemize}
Note that we always have $\tau_e\subseteq \tau(A)\subseteq \tau_s$ where $\tau_e$ is the euclidean topology and $\tau_s$ is the Sorgenfrey one.

We denote the space $(\mathbb{R}, \tau(A))$ as $H(A)$ and refer to it as the Hattori space associated to $A$.

It is well-known that the Euclidean and Sorgenfrey topologies have a somewhat curious relationship – while they share several topological properties, in others, they are entirely opposite. This prompts a natural question: under what conditions do Hattori's spaces preserve the properties of either of these topologies? Examples of this can be found in papers \cite{completez} and \cite{Hattori}, where two particularly notable patterns emerge: \begin{itemize}
    \item If the property is shared by both topologies, then any Hattori space maintains the property.
    \item If the complement of $A$ is countable, then the Hattori space associated with $A$ maintains the metric/completness properties of the euclidean topology (see \cite{completez}).
    
\end{itemize}

These two patterns raise some natural questions: \begin{itemize}
    \item Since both the Euclidean and Sorgenfrey topologies are quasimetrizable, is every Hattori space quasimetrizable? Note that any answer to this question breaks one of the previous patterns.
    \item Does these patterns hold true if we consider the space of continuous functions?
\end{itemize} 
In this paper we work on the last questions, more specifically, we try to determine when is $C_p(H(A))$  a normal/Lindelöf space, and along the way, we construct closed and discrete families in  $C_p(H(A))$, this could be considered the main focus of this paper. Also we find under which conditions $H(A)$ a quasimetrizable/semi-estratifiable space, and we give a characterization of compact sets in $H(A)$. We start by fixing some notation and recalling some definitions.\\

\section{Preliminaries}

Throughout this paper we will denote the Sorgenfrey line  as $\mathbb{S}$, meanwhile we will use $\mathbb{S}^\star$ to refer to the ``inverted" Sorgenfrey line, that is, the space that have as local bases intervals of the form $(a,b]$. As always, $\mathbb{R}$ represents the real line with his usual topology. $C_p(X)$ represents the space of real valued continuous functions defined on $X$ with the topology of pointwise convergence. Lastly, as we said before, $H(A)$ represents the Hattori space associated with $A$.  \\

We refer to the intervals $(x-\varepsilon,x+\varepsilon)$ and $[x,x+\varepsilon)$ as $B_e(x,\varepsilon)$ and $B_s(x,\varepsilon)$, respectively. Also, for a fixed $A\subseteq \mathbb{R}$ we define:\[B_A(x,\varepsilon)=\begin{cases}
    B_e(x,\varepsilon), \ \text{if} \ x\in A \\
    B_s(x,\varepsilon), \ \text{other case}
\end{cases}\]

The symbols $\text{cl}_e$, $\text{cl}_s$, $\text{cl}_{H(A)}$ will denote the closures in the euclidean topology, the Sorgenfrey one, and in $H(A)$, respectively.

We will need the concepts of quasimetrizable, $\gamma$-space, semi-stratifiable space, Moore space, $\sigma$-space, $\Sigma$-space, $p$-space, quasicomplete space and Čech complete space. We will use most of these to state exactly one theorem, so we prefer to refer the reader to \cite{handbook} (explicitly, chapter 10) for the first 6 definitions, \cite{creede} for the definitions of $p$-space and quasicomplete and lastly, \cite{engelking} for the definition of Čech complete space.

The main results concerning all of this definitions is summarized in the following proposition:
\begin{pro}\label{Relaciones}
    \begin{itemize}
        \item Every Moore space is both a quasimetrizable space and a semi-stratifiable space, and every quasimetrizable space is a quasimetrizable space. 
        \item Every Cech-Complete space is a $p$-space, and every $p$-space is quasicomplete.
        \item A Tychonoff space is a Moore space if and only if it is a semi-stratifiable $p$-space \cite{creede}.
        \item A space $X$ is a $\sigma$-space if and only if $X$ is a $\Sigma$-space with a $G_\delta$ diagonal.
    \end{itemize}
\end{pro}
Lastly an uncounable regular cardinal $\kappa$ is a caliber of a space $X$ if, for any family $\U\subseteq \tau\setminus \lbrace \emptyset\rbrace$ of cardinality $\kappa$, there exists $\U^\prime \subseteq \U$ such that $|\U^\prime|=\kappa$ and $\bigcap \U^\prime\neq \emptyset.$
All undefined terms will be interpreted as in \cite{engelking}.\\

\section{Metric-like properties of $H(A)$ and some compactness results}

We start our study of $H(A)$ by answering the following question: Under which conditions on $K\subseteq \mathbb{R}$, $K$ is a compact set in $H(A)$? Since $\tau(A)$ contains the euclidean topology, if $K$ is compact in $H(A)$ it is compact in $\mathbb{R}$, so we need $K$ to be a closed and bounded set in $\mathbb{R}$. We have the following full caracterization of compact sets in $H(A)$:

\begin{pro}
    
 Let $A, K \subseteq \mathbb{R}$ with $|K|\geq \aleph_0$. The following are equivalent:\begin{enumerate}
     \item $K$ is compact in $H(A)$.
 
\item  \begin{itemize}
    \item $K\setminus A$ is countable.
    \item For every $x\in K\setminus A$ there exists $\varepsilon>0$ such that $(x-\varepsilon,x)\cap K=\emptyset$.
\end{itemize}
\end{enumerate}

\end{pro}

\begin{pf}
    
Start by supposing that $K$ is compact in $H(A)$. Since $\tau_e\subseteq \tau(A)$, $K$ is compact in the euclidean topology. Now let's verify that $K$ satisfies (1). Note that $K$ is compact and submetrizable (again, because $\tau_e\subseteq \tau(A))$, so $K$ is a compact space with a $G_\delta$ diagonal, hence a metrizable space (this is a classic result by Sneider \cite{sneider}, see also \cite[Ch.~9, \S ~2 ]{handbook}). So $K\setminus A$, as a subspace of a metric space, it's also metrizable. On the other hand, the topology he inherits as a subspace of $H(A)$ is the one he inherits as a subspace of $\mathbb{S}$ (this is because $K\setminus A\subseteq \mathbb{R}\setminus A$, see 2.1 in \cite{Hattori}). Since the only metrizable subspaces of $\mathbb{S}$ are countable we conclude that $K\setminus A$  satisfies the first part of (2).\\

\noindent Let's move on to the second part of (2). Let $x\in \mathbb{R}\setminus A$ and suppose, for the sake of contradiction, that for every $\varepsilon>0$, $(x-\varepsilon,x)\cap K\neq \emptyset$. In particular, for each $n\in \N$ take $x_n\in (x-\varepsilon,x)\cap K$ and wlog suppose that $x_n<x_{n+1}$. Take the following open cover of $K$, $\U=\lbrace (-\infty,x_2)\rbrace \cup \lbrace (x_{n-1},x_{n+1}) \ | \ n\geq 2\rbrace \cup \lbrace [x,\infty)\rbrace $ (note that the last set is open because $x\notin A$). $\U$ lacks finite subcovers since every element of our sequence belongs to exactly one element of $\U$, thus contradicting the compactness of $K$. With this, we conclude the (1)$\rightarrow (2)$ implication.\\

\

\noindent Now let's prove the reverse implication, suppose that  $K$ is a compact set in the euclidean topology and satisfies (1). Let $\U\subseteq \tau(A)$ be a open cover of $K$ in $H(A)$. For each $x\in K$ let $U_x\in \U$ such that $x\in U_x$. We can take  $\varepsilon_x>0$ that satisfies the following conditions:\begin{itemize}
    \item If $x\in A$, then $(x-\varepsilon_x,x+\varepsilon_x)\subseteq U_x$.
    \item If $x\notin A$, then $[x,x+\varepsilon)\subseteq U_x$ and  $(x-\varepsilon,x)\cap K=\emptyset$. ($\star$)
\end{itemize}
Let's consider $\mathcal{V}=\lbrace B_e(x,\varepsilon_x) \ | \ x\in K\rbrace$. This is an open cover for $K$ consisting of open sets in the euclidean topology, and therefore exists a finite set $F\subseteq K$ such that \[
K\subseteq \bigcup_{x\in F} B_e(x,\varepsilon_x)
\]But ($\star$) implies that, if $x\in K\setminus A$,  then we can remove the left side of the interval $(x-\varepsilon_x,x+\varepsilon_x)$ without removing points of $K$, thus:
\[K\subseteq \bigcup_{x\in F} B_A(x,\varepsilon_x)\subseteq \bigcup_{x\in F} U_x\]
And so, we have found a finite subcover of $\U$. With this, we conclude that $K$ is compact in $H(A)$.
\end{pf} 

\noindent The first part of condition (1) does not come as a surprise because it is often encountered when demonstrating that specific properties of $\mathbb{R}$ are preserved in $H(A)$. On the other hand, the second part may appear sudden or unrelated. However, upon closer examination and with the following reinterpretation in mind, it does make sense:
 Condition (1) states that if $x\in K\setminus A$, then $x$ is not an accumulation point of $K$ in $Y$, and so we can restate our Proposition as follows:
\begin{pro}
    Let $A,K\subseteq \mathbb{R}$ with $|K|\geq \aleph_0$. Then $K$ is compact in  $H(A)$ if and only if $K$ is compact in the euclidean topology, $K\setminus A$ is countable and  $\text{der}_{Y}(K)\cap (K\setminus A)=\emptyset $.
\end{pro}

Since $H(A)$ is first-countable and hereditarily Lindelöf, we have actually proven a bit more, namely:
\begin{thm} Let $A,K\subseteq \mathbb{R}$ with $K$ infinite and compact in the euclidean topology. Then the following are equivalent:\begin{itemize}
    \item $K$ is compact in $H(A)$.
    \item $K$ is countably compact in  $H(A)$.
    \item $K$ is sequentially compact in $H(A)$.
    \item $|K\setminus A|\leq \aleph_0$ and $\text{der}_Y(K)\cap (K\setminus A)=\emptyset$.
\end{itemize}
\end{thm}

Now we move on to the question: When is $H(A)$ a quasimetrizable/semi-stratifiable? The first question will be answered with a classical condition, namely, if $|\mathbb{R}\setminus A|\leq \aleph_0$. The second question will require a  more elaborate condition. Let's start by addressing the semi-stratifiable case with a couple of simple results:

\begin{lem}
    Let $A\subseteq \mathbb{S}$ with a countable network. Then $A$ itself is countable.
\end{lem}
Take a subset $A$ of $\mathbb{S}$ and suppose that $A$ is semi-stratifiable. Note that $A$ is a quasimetrizable space (since $\mathbb{S}$ is), so $A$ is a semi-stratifiable, regular and quasimetrizable space, thus is a Moore space (See \cite[Ch.~10, \S~8]{handbook}). Since $\mathbb{S}$ is hereditarily Lindelöf, $A$ is a Lindelöf Moore space, thus second countable and the last Lemma allow us to conclude:

\begin{pro}\label{semies}
    Let $A\subseteq \mathbb{S}$ be semi-stratifiable. Then $A$ is countable. 
\end{pro}

The following Theorem can be obtained by using the results of \cite[\S~3]{completez} and Proposition \ref{Relaciones}. It contains most of the results concerning metric/completeness properties of Hattori spaces. 

\begin{thm}
    Let $A\subseteq \mathbb{R}$. The following are equivalent:\begin{itemize}
        \item $\mathbb{R}\setminus A$ is countable.
        \item $H(A)$ is Polish.
        \item $H(A)$ is completely metrizable.
        \item $H(A)$ is Čech-Complete.
        \item $H(A)$ is a Moore space.
        \item $H(A)$ is a p-space.
        \item $H(A)$ is quasicomplete.
        \item $H(A)$ is a $\sigma$-space.
        \item $H(A)$ is a $\Sigma$-space.
    \end{itemize}
\end{thm}

We are ready to add ``is semi-stratifiable" to these equivalences. Note that if $H(A)$ is a semi-stratifiable space, then $\mathbb{R}\setminus A$ is semi-stratifiable too, but this set inherits the same topology as a subspace of $\mathbb{S}$, thus Proposition \ref{semies} allow us to conclude that $\mathbb{R}\setminus A$ is countable. The other implication is trivial, so we get:

\begin{pro}
    Let $A\subseteq \mathbb{R}$. Then $H(A)$ is semi-stratifiable iff $\mathbb{R}\setminus A$ is countable.
\end{pro}

\noindent Since both $\mathbb{R}$ and $\mathbb{S}$ are quasimetrizable spaces, intuition suggests that $H(A)$ should be quasimetrizable regardless of the choice of $A$. This is similar to other properties shared between $\mathbb{R}$ and $\mathbb{S}$, such as being hereditarily Lindelöf and separable, being a Baire space, etc. However, this is not the case here, and the proof of this fact can be found in \cite{bennettquasi} and \cite{Quasi}, where Bennett proves what is known as the ``$\gamma$-space conjecture" (see Theorem 3.1 of \cite{bennettquasi}) for separable generalized ordered spaces and Kolner later modified said Theorem a little (see Theorem 10 of \cite{Quasi}). We must mention that the credit for following these ideas belongs to Li and Lin in \cite{lin}, and we will expand on their ideas by adding a couple of simple results connecting Benett and Kolner results. We present Bennett's theorem adapted to our context.

\begin{thm}
    Let $A\subseteq \mathbb{R}$. The following are equivalent:
    \begin{itemize}
        \item[a)] $H(A)$ us a $\gamma$-space.
        \item[b)] There exists a sequence of sets $(R_n)_{n=1}^\infty$ such that:\begin{itemize}
            \item $\mathbb{R}\setminus A=\bigcup_{n=1}^\infty R_n$
          \item  For each $p\in A\cap \text{cl}_e(R_n)$, there exists $y<p$ such that \\ $(y,p)\cap \text{cl}_e(R_n)=\emptyset$.\end{itemize}
        \item[c)] $H(A)$ is quasimetrizable.
    \end{itemize}
\end{thm}

This technically completely solves our question, but it is not very enlightening. Upon closer analysis of Bennett's proof, it turns out that we can deduce more, leading to the following result, of which we will prove the first implication for the sake of completness:

\begin{thm}\label{Quasi}
    Let $A\subseteq \mathbb{R}$. The following are equivalent:\begin{itemize}
        \item[a)]  $H(A)$ is a $\gamma$-space.
        \item[b)] There exists a sequence of sets $(R_n)_{n=1}^\infty$ such that:\begin{itemize}
            \item $\mathbb{R}\setminus A=\bigcup_{n=1}^\infty R_n$
          \item  For each $p\in A\cap \text{cl}_e(R_n)$, there exists $y<p$ such that \\ $(y,p)\cap \text{cl}_e(R_n)=\emptyset$.
            \item If $(x_n)\subseteq R_m$ is an increasing sequence, then it does not converge in $H(A)$.
       \end{itemize}
        \item[c)] $H(A)$ is quasimetrizable.
    \end{itemize}
\end{thm}

\begin{pf} Let $g:\N\times H(A) \rightarrow \tau_A $ be a $g$ function that satisfies the conditions of Definition \ref{Gamma}, let $Q=\lbrace q_k \ | \ k \in \N\rbrace$ be an enumeration of the rationals and  $Q_k=\lbrace q_i \ | \ i\leq k\rbrace$.  For every $n,m,k\in\N$ we define: \[
R(n,m,k)=\lbrace x\in \mathbb{R} \ | \ g(n,x)\subseteq [x,\infty) \ \& \ \alpha(n,x)=m \ \& \ g(m,x)\cap Q_k\setminus \lbrace x\rbrace\neq \emptyset\rbrace
\]

Let's see that $\mathbb{R}\setminus A=\bigcup_{(n,m,k)\in \N^3} R(n,m,k)$. Given $x\in \mathbb{R}\setminus A$, we know that $[x,\infty)$ is an open set, and since $g$ is a $\gamma$-function, this implies that there exists $n\in\N$ such that $g(x,n)\subseteq [x,\infty)$. This gives us our first parameter. For the second parameter, we can take $m=\alpha(n,x)$. Finally, since $\Q$ is a dense set in $H(A)$ and $g(m,x)\setminus \lbrace x\rbrace$ is a non-empty open set, there exists $k\in \N$ such that $Q_k\cap g(m,x)\setminus \lbrace x\rbrace\neq \emptyset$. Thus, $x\in R(n,m,k)$. \\

\noindent The other contention is simple. This is because the first condition defining $R(n,m,k)$ implies that $[x,\infty)$ is an open set, and therefore $x\in \mathbb{R}\setminus A$.\\

\noindent Let's consider $(x_j)\subseteq R(n,m,k)$, a sequence such that $x_j<p$ for all $j\in \mathbb{N}$. We will see that $(x_j)$ cannot converge to $p$ in $H(A)$. If $p\in \mathbb{R}\setminus A$, then it is clear that $(x_j)$ cannot converge to $p$ in $H(A)$. Therefore, necessarily $p\in A$, let's suppose that the sequence does converge.  \\

\noindent \textbf{Claim:} 
One of the following cases applies:\begin{itemize}
    \item There exists $q\in Q_k$ such that $p\in (-\infty,q]$ and $(-\infty,q] \cap Q_k\setminus \lbrace p\rbrace=\emptyset$.
    \item There exist $a,q\in Q_k$ such that $b<q$, $p\in (b,q]$ and $(b,q]\cap Q_k\setminus \lbrace p\rbrace=\emptyset$.
\end{itemize}

\noindent First, let's see that it is impossible that for every $q\in Q_k$, $q<p$. 
Let's assume that it happens, wlog $q_k=\max Q_k$. Since $(x_j)$ converges to $p$, there exists $N\in \N$ such that, $x_s\in (q_k,p)$ for every $s\geq N$. Since $x_s\in R(n,m,k)$ we conclude  $g(n,x_s)\subseteq [x_s,\infty)$, $\alpha(n,x_s)=m$ y $Q_k\cap g(m,x_s)\setminus \lbrace x_s\rbrace\neq \emptyset$, take $q$ in the last intersection, it follows that $g(m,q_i)\subseteq g(n,x_s)\subseteq [x_s,\infty)$ and thus $q_i\in [x_s,\infty)$ which is a contradiction. We conclude that there exists some $q\in Q_k$ such that $p\leq q$. The cases from the claim come from the following: \begin{itemize}
    \item If for every $q\in Q_k$ we have $p\leq q$, it is enough to take $q=\min Q_k$ so the first case from the claim is satisfied.
    \item There exist $q_i,q_j\in Q_k$ such that $q_i<p\leq q_j$. In this case we take $a=\max \lbrace x\in Q_k \ | \ x<p\rbrace$ y $q=\min \lbrace x\in Q_k \ | \ p\leq x\rbrace$ and it follows that these points satisfy the second point from the claim.
\end{itemize}
 Continuing the proof, now we assert that there exist $N_1$ and $N_2$ in $\mathbb{N}$ such that:\begin{itemize}
     \item For all $i\geq N_1$, $q\in g(m,x_i)$ (this $q$ is the same from the corresponding case of the claim).
     \item For all $i\geq N_2$, $x_i\in g(m,p)$.
 \end{itemize}
We prove the existence of $N_1$ following the cases from the claim:
\begin{itemize}
    \item In this case $N_1=1$. For all $i\in \N$ we have  $g(m,x_i)\subseteq g(n,x_i)\subseteq [x_i,\infty)$, the first inclusion follows from $m=\alpha(n,x_i)$. Since $(-\infty,q)\cap Q_k\setminus \lbrace p\rbrace=\emptyset$ and $g(m,x_i)\cap Q_k\setminus \lbrace x_i\rbrace\neq \emptyset$ we conclude that there exists $q_j\in Q_k\cap g(m,x_i)$ such that $q_j\geq q\geq p >x_i$ thus the interval $[x_i,q_j]$ is contained in $g(m,x_i)$, this because this last set is convex.  This guarantees $q\in g(m,x_i)$.
    \item By convergence, there exists $N_1\in \N$ such that, for all $i\geq N_1$ we have $x_i \in (a,p)\subseteq (a,q)$. Let $z\in Q_k\cap g(m,x_i)\setminus \lbrace x_i\rbrace$, since $g(m,x_i)\subseteq [x_i,\infty)$,  the election of $a$ and $q$ guarantees $z\geq q$. Thus $[x_i,z]\subseteq g(m,x)$ and it follows that $q\in g(m,x_i)$.
\end{itemize}
\noindent With this, we have proven the existence of $N_1$ for both cases. As for $N_2$, it is much simpler as we can directly apply the definition of convergence to the open set $g(m, p)$.\\

Now, let $i,j>N_1+N_2$ such that $x_i<x_j$ (we can always do this fixing $i$ and using convergence). Notice that $x_j,q\in g(m,x_j)$ and thus $[x_j,q]\subseteq g(m,x_j)$, from where we conclude that $p\in g(m,x_j)$. Since $g$ is a $\gamma$-function, we deduce \[g(m,p)\subseteq g(n,x_j)\subseteq [x_j,\infty)\]
Let's remember that $x_i\in g(m,p)$, thus $x_j\leq x_i$, which is a contradiction. Therefore, $(x_j)$ does not converge in the euclidean sense to $p$, but since $p\in A$, this implies that the sequence does not converge to $p$ in $H(A)$. \\

\noindent With what we have done so far, we have concluded that if $(x_n)$ is a sequence contained in $R(n,m,k)$, and $x_i < p$ for all $i$, then the sequence does not converge to $p$ in $H(A)$. As a particular case, no increasing sequence contained in $R(n,m,k)$ can converge in the Euclidean sense. \\

\noindent  Moreover, if $p\in A\cap \text{cl}_e(R(n,m,k))$, then there exist $\varepsilon>0$ such that $(p-\varepsilon,p)\cap R(n,m,k)=\emptyset$, on the contrary we could construct a sequence fully contained in $R(n,m,k)$ that converges to $p$ and  stays strictly to the left of $p$, but we just proved that this is imposible. notice that for this $\varepsilon>0$ we have $(p-\varepsilon,p)\cap \text{cl}_e(R(n,m,k))=\emptyset$. With this, we have simultaneously proven that both conditions we were looking for are satisfied for each $R(n,m,k)$, so we conclude this proof.
\end{pf}

\noindent This new version of the Theorem tells us a bit more, but it is still not clear enough which sets satisfy b). However, upon closer analysis, we obtain the following:

\begin{pro}
    
Suppose that $H(A)$ is a quasimetrizable space, and let \\ $\mathbb{R}\setminus A=\bigcup_{n=1}^\infty R_n$ where each $R_n$ satisfies $b)$ from Theorem 11. Then  for every $n\in\N$, $R_n$ is a closed set in $\mathbb{S}^\star$. That is, $\mathbb{R}\setminus A$ is a $F_\sigma$ set  in $Y$.\end{pro}
\begin{pf}
    Let's assume that there exists $p\in \text{cl}_Y(R_n)\setminus R_n$. With this, for each $m\in \mathbb{N}$, we can choose a point $x_m\in (p-\frac{1}{m},p)\cap R_n$. Moreover, we can select these points in such a way that they  give us an increasing sequence contained in $R_n$ that converges to $p$, which is impossible. We conclude that such a $p$ does not exist, and therefore, $R_n$ is closed in $Y$.
\end{pf}
\noindent 
Now we have an important condition: the set $\mathbb{R}\setminus A$ must be an $F_\sigma$ set in $Y$. Notice also that if $\mathbb{R}\setminus A=\bigcup_{n=1}^\infty R_n$ is an $F_\sigma$ set in $Y$, then each $R_n$ satisfies part b) of Theorem 10. To see this, suppose $p\in A\cap \text{cl}_e(R_n)$. Then, $p\notin R_n=\text{cl}_Y(R_n)$, so there exists $y<p$ such that $(y,p]\cap R_n=\emptyset$, and thus $(y,p)\cap R_n=\emptyset$, which implies $(y,p)\cap \text{cl}_e(R_n)=\emptyset$. With all this in mind, we can now give a definitive answer to our original question. Note that the following is essentially Theorem 10 of \cite{Quasi}.
\begin{thm}
    Let $A\subseteq \mathbb{R}$. The following are equivalent:\begin{itemize}
        \item $H(A)$ is a $\gamma$-space.
        \item $\mathbb{R}\setminus A$ is a $F_\sigma$ in $Y$.
        \item $H(A)$ is a quasimetrizable space.
    \end{itemize}
\end{thm}

\noindent To conclude this section, let's consider an explicit example of a set $A \subseteq \mathbb{R}$ such that $H(A)$ is not quasimetrizable. According to our theorem, it suffices to find a set $B$ that is not an $F_\sigma$ set in $\mathbb{S}^\star$, which is a simpler task than the original problem.
\begin{exa}A Hattori space that is not quasimetrizable. \\
\noindent Note that since $\mathbb{S}$ is a Baire space, $\mathbb{S}^\star$ is also a Baire space. Moreover, $\mathbb{R}\setminus \Q$ is a $G_\delta$ set in $Y$ because it is a $G_\delta$ set in the Euclidean topology. From these two observations, we can conclude that if we simply follow the usual proof that $\Q$ is not $G_\delta$ in $\mathbb{R}$, we can deduce that $\Q$ is not $G_\delta$ in $\mathbb{S}^\star$, and therefore $\mathbb{R}\setminus \Q$ is not an $F_\sigma$ set in $\mathbb{S}^\star$. Hence, we conclude that $H(\Q)$ is not a quasimetrizable space.
    
\end{exa}

\section{Lindelöf property and normality of $C_p(H(A))$}

Let's move on to the question of when $C_p(H(A))$ is a Lindelöf space. In order to prove that $C_p(H(A))$ is a Lindelöf space, it is enough to prove that it has countable network weight. Our next Proposition gives us a clear answer to this problem.
\begin{pro}\label{peso red}
    Let $A\subseteq \mathbb{R}$. Then $\text{nw}(H(A))=|\mathbb{R}\setminus A|+\aleph_0$.
\end{pro}
\begin{pf}
    One inequality is trivial, so let's prove the other one. Let $\mathcal{B}$ be a network for $H(A)$, let's see that $|\mathcal{B}|\geq |\mathbb{R}\setminus A|$. For each $x\in \mathbb{R}\setminus A$ we have $[x,x+1)\in \tau(A)$ and therefore it exists $B_x\in \mathcal{B}$ such that $x\in B_x\subseteq [x,x+1)$, this is a one to one assignation, so $|\mathcal{B}|\geq \mathbb{R}\setminus A$. Thus $\text{nw}(H(A))\geq \mathbb{R}\setminus A+\aleph_0$.
\end{pf}

Thanks to the previous proposition and recalling that $\text{nw}(X)=\text{nw}(C_p(X))$, we obtain:
\begin{thm}
    Let $A\subseteq \mathbb{R}$.  If $|\mathbb{R}\setminus A|\leq \aleph_0$ then $C_p(H(A))$ is Lindelöf.
\end{thm}

Moreover, we can do more. First, let's recall the following result (see \cite{tkachuk2011cp}, problem 249):

\begin{pro}Let $\omega_1$ be a caliber for a space $X$. Then $C_p(X)$ is a Lindelöf $\Sigma-$space if and only if $X$ has a countable network.\end{pro}
    Since $H(A)$ is separable, $\omega_1$ is a caliber for $H(A)$, thus we can use the last result and Proposition  \ref{peso red} to conclude:

\begin{thm}
    Let $A\subseteq \mathbb{R}$. Then $C_p(H(A))$ is a Lindelöf $\Sigma-$space if and only if $\mathbb{R}\setminus A$ is countable.
\end{thm}

Let's move on with the normality of $C_p(H(A))$. Since the condition  $ |\mathbb{R} \setminus A|\leq \aleph_0$ seems to preserve many properties of $\mathbb{R}$ to $H(A)$, the natural question would be:
\begin{question} Let $A\subseteq\mathbb{R}$. Is $C_p(H(A))$ a normal space if and only if $|\mathbb{R}\setminus A|\leq \aleph_0$?
\end{question}

In the following we will try to solve this question. One implication is  straightforward, if $|\mathbb{R}\setminus A|\leq \aleph_0$, then $H(A)$ is a separable metric space and thus $C_p(H(A))$ is normal. So, the real question is:

\begin{question}
    Let $A\subseteq \mathbb{R}$. If $C_p(H(A))$ is normal, ¿$|\mathbb{R}\setminus A|\leq \aleph_0$?
\end{question}

\noindent To tackle this problem, we will proceed by contrapositive. If the set $\mathbb{R}\setminus A$ satisfies certain conditions (in addition to being uncountable), then $C_p(H(A))$ contains an uncountable closed discrete set, and therefore $C_p(H(A))$ is not normal by Jones Lemma.\\

\noindent One way to prove that $C_p(\mathbb{S})$ is not normal is by showing that it contains a closed and discrete subspace of size $\mathfrak{c}$. Explicitly, for every $a\in (0,1)$ define the following function:  \[
f_a(x)=\begin{cases} 0 \ \ \ & x\in(-\infty,-1) \\
1  &x\in[-1,-a)\\
0  &x \in [-a,a)\\
1  &x\in [a,1)\\
0  &x\in [1,\infty)
\end{cases}
\]
In this way, $\lbrace f_a \ | \ x\in(0,1)\rbrace$ is a closed and discrete subspace of $C_p(\mathbb{S})$ with cardinality $\mathfrak{c}$ (a proof of this fact can be found in \cite{BLOG})\\

We can replicate the previous construction on any interval which have the sorgenfrey topology, so we get the following partial results:

\begin{pro} Let $A\subseteq \mathbb{R}$. If there exists $a<b$ such that $[a,b]\subseteq \mathbb{R}\setminus A$, then $C_p(H(A))$ contains a closed and discrete subspace of cardinality $\mathfrak{c}$, and thus is not normal space.
    
\end{pro}
\noindent Containing a non trivial closed interval is equivalent to containing a non trivial open interval, so:
\begin{pro}
    Let $A\subseteq \mathbb{R}$. If $\text{int}_e(\mathbb{R}\setminus A)\neq \emptyset$, then $C_p(H(A))$ contains a closed and discrete subspace of cardinality $\mathfrak{c}$, and thus is not normal space.
\end{pro}

 \noindent Particularly, if $A$ is a closed set (and it is not $\mathbb{R}$) then $C_p(H(A))$ is not a normal space. \\ 
From here on, the idea is to gradually weaken the condition of containing intervals as much as possible to reach a more general result. Let's start with the following:

 \begin{pro}\label{Simetricos densos}
     Let $A\subseteq \mathbb{R}$ such that  $B=\mathbb{R}\setminus A$ is  symmetric, meaning that if $x\in B$, then $-x\in B$, and $|B|=\mathfrak{c}$ . Additionally, let's assume that $B\cap (0,1)$ is dense in $(0,1)$ and there exists $q\in B\cap [1,\infty)$. Then $C_p(H(A))$ contains a closed and discrete subset of size $\mathfrak{c}$, and therefore it is not a normal space.
 \end{pro}
\begin{pf}
    Let $B=\mathbb{R}\setminus A$. For each $a\in (0,1)\cap B$ we define:  

\[
f_a(x)=\begin{cases} 0 \ \ \ & x\in(-\infty,-q) \\
1  &x\in[-q,-a)\\
0  &x \in [-a,a)\\
1  &x\in [a,q)\\
0  &x\in [q,\infty)
\end{cases}
\]

\noindent By defining $f_a$ in this way, we note that the ``jumps" occur at points where the local topology is that of the Sorgenfrey Line, so we do not break continuity in $H(A)$, so $\mathcal{F}\subseteq C_p(H(A))$. Let's proceed to show that it is indeed a closed and discrete set.\\

\noindent Let $a\in B\cap (0,1)$. Let's consider the basic open set $V=\left[f_a,\lbrace a,-a\rbrace, \frac{1}{2}\right]$ and let's verufy the following cases;\\

\noindent\textbf{Case 1.-} $0<b<a<1$. Here we have that $f_b(-a)=1$, while $f_a(-a)=0$, so $f_b\notin  V$.\\

\noindent\textbf{Caso 2.-} $0<a<b<1$. This time $f_b(a)=0$ but $f_a(a)=1$, and so, once again, we conclude that $f_b\notin V$.\\

\noindent From the previous cases we conclude that $V\cap \mathcal{F}=\lbrace f_a\rbrace$, so $\mathcal{F}$ is a discrete subset of $C_p(H(A))$. Proving that it is closed is the difficult part. Note that our functions only take values in  $\lbrace 0,1\rbrace$, if $g\in C_p(H(A))$ is such that $g(z)\notin \lbrace 0,1\rbrace$ for some  $z\in \mathbb{R}$, then $U=[g,\lbrace z\rbrace,\min \lbrace |g(z)|,|1-g(z)|\rbrace]$ satisfies $U\cap \mathcal{F}=\emptyset$. With the previous argument, we can suppose wlog that $g$ only takes the values $0$ or $1$. Once again, we have two cases:\\

\noindent \textbf{Case 1.-} There exists  $a\in B\cap (0,1)$ such that $g(a)\neq g(-a)$. Here we have two sub cases:\\

\noindent \textbf{Sub case 1.1.-} $g(a)=0$ y $g(-a)=1$. Let's take $V=[g,\lbrace \pm a\rbrace, \frac{1}{2}]$. Suppose  there exists $b\in B\cap(0,1)$ such that $f_b\in V$. Then we have that $f_b(a)\in (-\frac{1}{2},\frac{1}{2})$, that is, $f_b(a)=0$, analogously, $f_b(-a)=1$. $f_b(a)=0$ and $a>0$ imply that $a\in (0,b)$, so $a<b$. On the other hand $f_b(-a)=1$ implies that $-a\in[-q,-b)$, so $-a<-b$ and thus $a>b$. Since these two inequalities contradict each other, we deduce that such $b$ cannot exist, so  $V\cap \mathcal{F}=\emptyset$.\\

\noindent \textbf{Sub case  1.2.-} $g(a)=1$ y $g(-a)=0$. Once again let's take $V=[g,\lbrace \pm a \rbrace, \frac{1}{2}]$ and suppose  there exists $b\in B\cap (0,1)$ such that $f_b\in V$. Similar reasoning to the previous paragraph allows us to deduce that $f_b(a)=1$ and $f_b(-a)=0$, so $a\in [b,q)$ and $-a\in [-b,b)$, that is, $a\leq b$ and $-b\leq a$, thus $a=b$. With this we obtain $V\cap \mathcal{F}\subseteq \lbrace f_a \rbrace$, since $C_p(H(A))$ is Hausdorff we can find an open set $U$ such that $g\in U$ but $f_b\notin U$, so $V\cap U\cap \mathcal{F}=\emptyset$ and we conclude this sub case.\\

\noindent \textbf{Case 2.-} For all $y\in B\cap(0,1)$, $g(y)=g(-y)$. Once more, we have sub cases:\\

\noindent \textbf{$g$ is constant in $(0,1)$}\\

\noindent \textbf{Sub case 2.1.1-} $g(x)=0$ for all $x\in (0,1)$.\\
By continuity we obtain $g(-1)=0$, now  let's define $V=[g,\lbrace -1\rbrace,\frac{1}{2}]$. For every $b\in B\cap (0,1)$, $f_b(-1)=1$ so $f_b\notin V$ and we conclude.\\

\noindent \textbf{Sub case 2.1.2-} $g(x)=1$ for all $x\in (0,1)$.\\
It suffices to take $V=[g,\lbrace 0\rbrace, \frac{1}{2}]$ and reason in the same way as the last sub case noting that $f_b(0)=0$ for all $b\in B\cap (0,1)$.\\

\noindent \textbf{$g$ is not constant in $(0,1)$}\\

\noindent \textbf{Claim:} There exist $c,d\in B\cap (0,1)$ such that $g(c)\neq g(d)$.\\
Since $g$ is not constant in $(0,1)$ there exist $x,y\in (0,1)$ such that $g(x)\neq g(y)$, due to continuity we can find $\varepsilon>0$ cuch that $g$ is constant in $[x,x+\varepsilon)$ and $[y,y+\varepsilon)$ (it does not matter if the points have the euclidean or Sorgenfrey topology). Thanks to density we can find $c\in [x,x+\varepsilon)\cap B\cap (0,1)$ and $d\in[y,y+\varepsilon)\cap B\cap (0,1)$, it follows that $g(c)\neq g(d)$. Moreover, wlog $c<d$.\\ 

\noindent \textbf{Sub case 2.2.1} $g(c)=0$ and $g(d)=1$. \\
Let's define $\U=\lbrace a \in B\cap (0,d) \ | \ g(a)=0\rbrace$ and let $u= \sup \U$. Notice that  $u\leq d$ and let's see that $u\in B$. Suppose that $u\notin B$, since $g$ is continuous there exists $\varepsilon>0$ such that $g[(u-\varepsilon,u+\varepsilon)]\subseteq \lbrace g(u)\rbrace$, we take $z\in (u-\varepsilon,u)\cap \U$, then $0=g(z)=g(u)$. By density of $B$ we find an $w\in B\cap (0,1)\cap (u,u+\varepsilon)$ and thus $g(w)=0$, a contradiction. We conclude that $u\in B$.\\

\noindent Moreover, $g(u)=1$, since in other case we could find some point in $(u,d)\cap B\cap (0,1)$ such that $g(u)=0$ contradicting the choice of $u$ (if $u=d$ there was nothing to do). Construct an increasing sequence $(x_n)\subseteq U\cap (0,1)$ converging to $u$, it follows that $g(x_n)=0$ for all $n$, then $(-x_n)$ is a decrecing sequence that converges to $-u$ and  $g(x_n)=g(-x_n)=0$, due to continuity we conclude $g(-u)=0$, but, since $g(u)=1$, we also have $g(-u)=1$, a contradiction. We deduce that this case is impossible.\\

\noindent \textbf{Sub case 2.2.2} $g(c)=1$ y $g(d)=0$. \\
Analogously to the previous case, we define $\U=\lbrace a \in B\cap (0,d) \ | \ g(a)=1\rbrace$ and take $u=\sup \U$, once again $u\leq d$. Suppose $u\notin B$, by continuity there exists $\varepsilon>0$ such that $g[(u-\varepsilon,u+\varepsilon)]\subseteq \lbrace g(u)\rbrace$. Since there exists $z\in U\cap (u-\varepsilon,u)$ we conclude that $g(u)=1$, now we take  $w\in B\cap (0,1)\cap (u,u+\varepsilon)$, but this implies $g(u)=1$(note that $u+\varepsilon<d$, otherwise $g(d)=1$), a contradiction. So we conclude $u\in B$. From here onward evrything proceeds analogously to the previous case to get a contradiction.\\ 

\noindent Since the last two sub cases lead us to contradictions, we conclude that the condition $g(a)=g(-a)$ for all $a\in B$, implies that $g$ is constant in $(0,1)$, a case that we solved previously.

\end{pf}

\noindent The crucial parts for the construction were as follows, let $B=\mathbb{R}\setminus A$:\begin{itemize}
    \item $B$ is symmetric.
    \item $B$ is dense in some interval, say $(0,b)$.
    \item There exists $q\in [b,\infty)\cap B$.
\end{itemize}
\noindent Reading carefully we note that the second condition implies the third in the following sense: if $B\cap(0,b)$ is dense in $(0,b)$, then $B\cap (0,y)$ is dense in $(0,y)$ for every $y<b$ and  we know there exists $q\in (y,b)\cap B$ since $(y,b)$ is an open set of $(0,b)$.\\

\noindent We can generalize the result by extending the possible symmetries. Remember that the reflection with respect to a point $a \in \mathbb{R}$ is given by $r_a(x) = 2a - x$, and all these reflections satisfy the following property: they transform increasing sequences ``on the right of a" into decreasing sequences ``on the left of  a". This is a key argument for the final part of the previous proof. \\

We generalize the symmetry condition with the following definition:
\begin{definition} Let $A \subseteq \mathbb{R}$ and $z \in \mathbb{R}$. We say that $A$ is symmetric with respect to $z$ if for every $x \in A$, $r_z(x) \in A$.

\end{definition}
Also let's recall the following classic results (\cite{libroCp}):

\begin{pro}
    If $C_p(X)$ is normal, then it is collection-wise normal. 
\end{pro}

\begin{pro} Let $X$ collection-wise normal and ccc. Then $e(X)=\aleph_0$. 
    
\end{pro}
    
 The preceding results allow us to improve on our reasoning, we do not need a closed and discrete set of cardinality $\mathfrak{c}$, we only need that it is uncountable.\\ 

\noindent With this new ideas we can improve  Proposition 6 in the following way:

\begin{thm}
    Let $A\subseteq \mathbb{R}$ y $B=\mathbb{R}\setminus A$, suppose that there exists $(z-a,z+a)=I$ such that :\begin{itemize}
      \item $B\cap I$ is dense in $I$.
      \item $B\cap I$ is symmetric with respect to $z$.
      \item $|B\cap I|>\aleph_0$
    \end{itemize}
    Then $C_p(H(A))$ contains a closed and discrete subspace of size $|B\cap I|$ and thus is not normal.
\end{thm}

While the conditions of the theorem may seem difficult to achieve, it is relatively simple to construct sets $B$ that satisfy them. For example, we start with $C = \mathbb{Q} \cap (-1,1) \setminus {0}$. This set is already dense in $(-1,1)$ and symmetric with respect to zero. We can then add as many irrationals as we want while maintaining symmetry, and the result will still be dense. In particular, we have the following:
\begin{pro}

For every $\aleph_0 < \kappa \leq \mathfrak{c}$, there exists $B \subseteq \mathbb{R}$ that satisfies the conditions of the previous theorem and has a cardinality of $\kappa$. Moreover, there exist $\mathfrak{c}$ distinct sets that satisfy these conditions.
\end{pro}
\begin{pf}
    It suffices to perform the same construction while varying the symmetry point. In other words, we start with $\mathbb{Q} \cap (z-1,z+1) \setminus {z}$ and then add the desired irrationals. These sets will be distinct because the initial set of rationals is different for each $z$.
\end{pf}
\noindent Theorem 3 allows us to construct closed and discrete sets by leveraging the density and symmetry of a sufficiently large set. Now, let's maintain the symmetry but shift to the opposite side by considering sets that are nowhere dense.

\noindent Let $\Delta$ be the Cantor set with its usual construction, that is, $\Delta = \bigcap_{n=1}^\infty C_n$ where each $C_n$ is the union of $2^n$ intervals of length $\frac{1}{3^n}$. With this in mind, we will say that $u \in (0,1)$ is a left (right) endpoint of $\Delta$ if there exists $n \in \mathbb{N}$ such that $u$ is the left (right) endpoint of one of the intervals that form $A_n$. Furthermore, we write $[0,1]\setminus\Delta = \bigcup_{n=1}^\infty I_n$, where the $I_n$ are the intervals discarded when constructing $\Delta$. We will need the following results concerning $\Delta$.

\begin{lem}
Let $u\in \Delta$. Then: \begin{itemize}
    \item $u$ is a left endpoint of $\Delta$ iff there exists $u^- \in \Delta$ such that $u^-<u$ and $(u^-,u)\cap \Delta=\emptyset$.
    \item $u$ is a right endpoint $\Delta $ iff there exists $u^+\in \Delta$ such that $u^+>u$ and $(u,u^+)\cap \Delta \neq \emptyset$.
    \item If $u$ is not a left endpoint, then there exists $(x_n)\subseteq \Delta$ strictly increasing sequence such that $x_n\rightarrow u$ (in the euclidean sense). Moreover, for every $\varepsilon>0$, $(u-\varepsilon,u)\cap\Delta \neq \emptyset$.
    \item If $u$ is not a right endpoint, then there exists $(x_n)\subseteq \Delta$ strictly decreasing sequence such that $x_n\rightarrow u$ (in the euclidean sense). Moreover, for every $\varepsilon>0$, $(u,u+\varepsilon)\cap\Delta \neq \emptyset$.
    \end{itemize}
\end{lem}
\noindent Explicitly, $u^+$ is the left endpoint of the next interval that forms $A_n$ (ordered in increasing order), and similarly for $u^-$.

\begin{thm}
    Let $A\subseteq \mathbb{R}$ such that $\Delta \subseteq \mathbb{R}\setminus A$. Then $C_p(H(A))$ contains a closed and discrete set of cardinality $\mathfrak{c}$ and thus is not normal.
\end{thm}

\begin{pf}
    
Let $B=[\frac{1}{2}, 1) \cap \Delta$ and $f(x) = 1 - x$, which is the reflection with respect to $\frac{1}{2}$. Note that $f[\Delta] = \Delta$. For each $b \in [\frac{1}{2}, 1) \cap \Delta$, we define a function as follows:
    
\[
g_b(x)=\begin{cases} 0 \ \ \ & x\in(-\infty,0) \\
1  &x\in[0,f(b))\\
0  &x \in [f(b),b)\\
1  &x\in [b,1)\\
0  &x\in [1,\infty)
\end{cases}
\]
\noindent Note that the ``jumps" of $g_b$ occur at points in $\Delta$, so continuity is preserved. Therefore, $\mathcal{F}=\lbrace g_b \ | \ b\in B\rbrace \subseteq C_p(H(A))$. Additionally, by construction, $g_b$ is constant on $I_n$ for all $n \in \mathbb{N}$. Let's see that $\mathcal{F}$ is closed and discrete in $C_p(H(A))$.\\

\noindent Let's see that $\mathcal{F}$ is discrete.\\
\noindent Given $b\in B=\Delta \cap [\frac{1}{2},1)$, we define $U=[g_b,\lbrace b,f(b)\rbrace,\frac{1}{2}]$. Let's take $c\in B$, we have the following two cases:\begin{itemize}
    \item \textbf{$c<b$.} Then $f(b)<f(c)$ and thus $g_c(f(b))=1$, but $g_b(f(b))=0$, which implies $g_c\notin U$.
    \item \textbf{$c>b$.} This time $g_b(b)=1$ while $g_c(b)=0$, once again, $g_c\notin \mathcal{F}$.
\end{itemize}
From this two cases we deduce that $U\cap\mathcal{F}=\lbrace g_b\rbrace$ and thus $\mathcal{F}$ is a discrete subspace.\\

\noindent Now we only need to show that $\mathcal{F}$ is closed. For this, we will consider the following cases, let $g\in C_p(H(A))$.\\

\noindent \textbf{Case 1.} $g[\mathbb{R}]$ is not contained in $ \lbrace 0,1\rbrace$.\\
It follows that there exists $x\in \mathbb{R}$ such that $g(x)\notin \lbrace 0,1\rbrace$, so we can define \[\varepsilon=\min \lbrace |g(x)|, |1-g(x)|\rbrace>0\]
under this conditions, let's take $U=[g,\lbrace x\rbrace,\varepsilon]$, it follows trivially that for every  $b\in B$, $|g_b(x)-g(x)|\geq \varepsilon$ and thus $g_b\notin U$.\\

\noindent Henceforth, we assume that $g[\mathbb{R}]\subseteq \lbrace 0,1\rbrace$.\\

\noindent \textbf{Case 2.-}$g$ is not constant in the intervals discarded during the construction of $\Delta$. \\
Then there exists $n \in \mathbb{N}$ such that $g$ is not constant on $I_n$, so there exist  $x,y\in I_n$ such that $g(x)=0$ and $g(y)=1$. Let's take $U=[g,\lbrace x,y\rbrace,\frac{1}{2}]$. For every $b\in B$, since $g_b$  is constant in $I_n$, one of the quantities $|g(x)-g_b(x)|$ or $|g(y)-g_b(y)|$ is exactly  $1$ and thus $g_b\notin U$.\\

\noindent From here onwards, we will also assume that $g$ is constant on each interval $I_n$.\\

\noindent \textbf{Case 3.-} There exists $a\in B$ such that $g(a)\neq g(f(a))$.\\
First let's suppose $g(a)=0$ and $g(f(a))=1$. Let $U=[g,\lbrace a,f(a)\rbrace, \frac{1}{2}]$ y $b\in B$, then:\begin{itemize}
    \item If $b<a$, we get $g_b(a)=1$.
    \item If $b>a$, then $f(b)<f(a)$, so $g_b(f(a))=0$.
\end{itemize}
Both cases lead us to $g_b\notin U$ and thus $\mathcal{F}\cap U\subseteq \lbrace g_a\rbrace$, since the space is Hausdorff this suffices. The case $g(a)=1$ and $g(f(a))=0$ is analogous.\\

\noindent We add to our list of assumptions $g(a)=g(f(a))$ for all $a\in B$.\\

\noindent\textbf{Case 4.-} $g$ is constant in $B$.\\
\begin{itemize}
    \item $g(x)=0$ for all $x\in B$. By continuity and the fact that $g(f(b))=g(b)$ we deduce that $g(0)=0$. We take $V=[g,\lbrace 0\rbrace,\frac{1}{2}]$, since for all $b\in B$, $g_b(0)=1$ it follows that $g_b\notin U$.
    \item $g(x)=1$ for all $x\in B$. In this case we take $V=[g,\lbrace \frac{1}{2}\rbrace,\frac{1}{2}]$, for all $b\in B$, $g_b(\frac{1}{2})=0$ thus $g_b\notin U$.
\end{itemize}
In both cases we conclude. Lastly,\\

\noindent \textbf{Case 5.-} $g$ is not constant $B$.\\
First let's suppose that there exist $a,b\in B$ with $b<a$ such that $g(b)=0$ and $g(a)=1$ (the other case is analogous). Let's define the following set \[W=\lbrace c\in \Delta \ | \ c<a \ \& \ g(c)=0\rbrace\]
$b\in a$ so $W$ is non empty, on the other hand, $a$ is an upper bound for $W$, thus $\sup W=u$ exists, $b\leq u\leq a$ and $u$ is in $\Delta\cap [\frac{1}{2},a]$ since this set is closed in the euclidean topology.\\

\noindent \textbf{Claim:} $u$ is not a left endpoint.\\
Suppose that it is a left endpoint and let's see that this implies $u \neq a$. If $u = a$ and $u$ is a left endpoint, then $u^-$ is an upper bound for $W$ (since $g(u) = g(a) = 1$), which contradicts the choice of $u$. Now, since $u$ is a left endpoint, we can take a strictly decreasing sequence $(x_n) \subseteq B$ that converges to $u$ in the usual sense. Moreover, without loss of generality, we can assume that $x_n \in (u, a)$ for all $n \in \mathbb{N}$, which in turn allows us to assume that $g(x_n) = 1$ for all $n$ (otherwise, it would contradict the choice of $u$). The continuity of $g$ then implies that $g(u) = 1$, which again leads us to $u^-$ being an upper bound for $W$. Thus, $u$ cannot be a left endpoint.\\

\noindent \textbf{Claim:} $g(u)=0$. \\
\noindent Suppose $g(u) = 1$. Since $u$ is not a left endpoint, we can take a strictly increasing sequence $(x_n) \subseteq \Delta \cap [\frac{1}{2},1)$ that converges to $u$ in the usual sense. Moreover, as $u = \sup W$ and $u \notin W$, we can choose the sequence in such a way that $g(x_n) = 0$ for all $n$. It follows that $f(x_n)$ is a strictly decreasing sequence converging to $f(u)$ and such that $g(f(x_n)) = g(x_n) = 0$ for all $n \in \mathbb{N}$. This implies that $g(f(u)) = 0$ by continuity, and therefore $0 = g(f(u)) = g(u) = 1$, which is impossible. We conclude that $g(u) = 0$, which in turn allows us to conclude that $u < a$.\\

\noindent \textbf{Claim:} $u$  is a right endpoint.\\
\noindent Suppose this does not happen. By continuity, there exists $\varepsilon > 0$ such that $g[[u,u+\varepsilon)]\subseteq \lbrace 0\rbrace$, since $u$ is not a right endpoint 
\[(u,u+\varepsilon)\cap (u,a)\cap \Delta\neq \emptyset\] which contradicts the fact that $u$ is the supremum of $W$.\\

\noindent \textbf{Claim:} 
There does not exist a strictly increasing sequence $(x_n) \subseteq \Delta$ converging to $u$ such that $g(x_n) = 1$ for all $n \in \mathbb{N}$. That is, there exists $\varepsilon>0$ such that $g[(u-\varepsilon,u)]\subseteq \lbrace 0\rbrace$.\\
Suppose that such a sequence exists. Then, $(f(x_n))$ is a strictly decreasing sequence converging to $f(u)$. By the continuity of $g$, we can conclude that $1 = g(f(u)) = g(u) = 0$, which is a contradiction. From this point forward in the proof, we will use this $\varepsilon$.\\

\noindent Let's take $y\in(u-\varepsilon,u)\cap \Delta$, this is posible since $u$ is not a left endpoint, moreover, since $u$ is a right endpoint, we can consider $u^+$, and furthermore, we have that $u^+\leq a$ and $g(u^+)=1$. IN order to finish the proof let's prove that $U=[g,\lbrace u,u^+,y\rbrace,\frac{1}{2}]$ satisfies \[U\cap \mathcal{F}\subseteq \lbrace g_u,g_{u^+}\rbrace\]
which allow us to conclude since $C_p(H(A))$ is Hausdorff.\\

\noindent Let $c\in \Delta\cap[\frac{1}{2},1)$. We have three final sub cases:\\
\begin{itemize}
    \item $c<y$. In this case $g_c(y)=1$ and $g(y)=0$.
    \item $c\in[y,u)$. Here $g_c(u)=1$ and $g(y)=0$.
    \item $b>u^+$. We have $g_c(u^+)=0$ but $g(u^+)=1$,
\end{itemize}
In any case we deduce that $g_c\notin U$, and so we can conclude.
\end{pf}
The previous theorem has quite strong conclusions. Let's start by recalling a couple of results. The first one is a classical result concerning analytic spaces, and the second one is a theorem by Sorgenfrey.
\begin{thm}
    
Let $B \subseteq \mathbb{R}$ be an analytic and uncountable set. Then $B$ contains a set homeomorphic to the Cantor set $\Delta$.
\end{thm}

\begin{thm}
    Let $C, D \subseteq \mathbb{R}$ be compact nowhere dense sets. Then there exists an order isomorphism $f: C \rightarrow D$ between them. In particular, any two Cantor sets in $\mathbb{R}$ are order isomorphic.
\end{thm}

\noindent 
Thanks to these two results, we have that if $\mathbb{R} \setminus A$ is an analytic set (for example, open, closed, etc.), then it contains a copy of the Cantor set, let's call it $C$. Moreover, we can provide an order isomorphism between $\Delta$ and $C$, which allows us to "transfer" our construction of the usual Cantor set to this case in the following way:

\

\noindent 
Let $g$ be the order isomorphism. Since $g\left(\frac{1}{3}\right) < g\left(\frac{2}{3}\right)$, we can take an intermediate point, let's say $c$. This point will act as a "division" point, similar to how $\frac{1}{2}$ did originally. Now, we define $D = C \cap [c, \infty) \setminus {g(1)}$ and $I = (-\infty, c] \cap C \setminus {g(0)}$. These two sets are disjoint, and we have $C = D \cup I \cup {g(0), g(1)}$. We also need the following auxiliary function $\phi: C \rightarrow C$ given by $\phi(g(x)) = g(1-x)$. Since $g$ is bijective, $\phi$ is well-defined and bijective. Moreover, the fact that $g$ preserves the order implies that $\phi[D] = I$. Notice that $\phi$ transforms increasing sequences in $D$ into decreasing sequences in $I$. This function acts as a ``pseudo-symmetry," taking on the role of $f(x) = 1-x$ in the original construction.\\
\

\noindent The only thing we haven't reinterpreted is the concept of left or right endpoint, but we can do that very simply. A left/right endpoint of $C$ is the image under $g$ of a corresponding type of endpoint in $\Delta$. For example, $g\left(\frac{2}{3}\right), g\left(\frac{2}{9}\right)$ are left endpoints of $C$. All the conclusions of Lemma 1 can be naturally adapted to this definition using the fact that $g$ preserves order and is a homeomorphism.\\
\

Finally, for each $d \in D$, we define 

\[
h_d(x)=\begin{cases} 0 \ \ \ & x\in(-\infty,g(0)) \\
1  &x\in[g(0),\phi(d))\\
0  &x \in [\phi(d),d)\\
1  &x\in [d,g(1))\\
0  &x\in [g(1),\infty)
\end{cases}
\]

\noindent Under this reinterpretation of the key concepts used in the proof of Theorem 4, we can replicate that proof for this more general case by making the corresponding replacements to conclude that $\lbrace h_d \ | \ d \in D\rbrace$ is a closed and discrete set in $C_p(H(A))$. Therefore, we can state the following result:

\begin{thm}
    Let $A\subseteq \mathbb{R}$ such that $\mathbb{R}\setminus A$ is uncountable and analytic. Then $C_p(H(A))$ is not normal.
\end{thm}

Unfortunately, the previous theorem is not an ``if and only if" statement, as the following example shows:

\begin{exa}Let $B$ be a Bernstein set that is also an additive subgroup of $\mathbb{R}$ (see \cite{bernstein} and \cite{bernstein2}). It follows that $B$ is an uncountable set that is symmetric with respect to $0$, dense in $\mathbb{R}$, but not analytic. Therefore, by Proposition \ref{Simetricos densos}, we can conclude that $C_p(H(\mathbb{R}\setminus B))$ is not normal, even though $B$ is not analytic.
\end{exa}

All these partial advances naturally lead to the following question:
\begin{question}
    If $\mathbb{R}\setminus A$ is uncountable, then, is ¿$C_p(H(A))$ not normal?
\end{question}

\noindent One of the most significant models of ZF is Solovay's model, in which the Axiom of Choice is not satisfied fully (only the Axiom of Dependent Choice, DC), and furthermore, in this model, $\mathbb{R}$ has the  Perfect Set Property, meaning that every infinite subset is either countable or contains a perfect set. Since every perfect set is analytic, this implies that every non-countable subset of $\mathbb{R}$ contains a copy of $\Delta$. Therefore, in this model, we have an affirmative answer to the previous question:

\begin{thm}
    (In Solovay Model) Let $A\subseteq \mathbb{R}$. Then if $\mathbb{R}\setminus A$ is uncountable then $C_p(H(A))$ is not normal.
\end{thm}

\noindent From the last Theorem it follows that, if there exist an uncountable $B$  such that  $C_p(H(B))$ is normal, then it would need AC for his construction. Thus it will not be easy to describe such a set.

\bibliographystyle{abbrv}
\bibliography{biblio}
\end{document}